\documentclass[11pt]{article}
\usepackage{geometry, url}

\input epsf.sty
\epsfverbosetrue 

\usepackage{latexsym}
\usepackage{amsfonts}
\usepackage{amssymb}
\usepackage{amsmath}

\newcommand{\R}{\mathbb{R}}
\newcommand{\Q}{\mathbb{Q}}
\newcommand{\Z}{\mathbb{Z}}
\newcommand{\C}{\mathbb{C}}
\newcommand{\IP}{\mathbb{P}}

\newcommand{\N}{\mathbb{N}}

\newcommand{\rs}{\mbox{$\widehat{\C}$}}

\def\SSS{{\mathcal S}}
\def\TTT{{\mathcal T}}
\def\AAA{{\mathcal A}}

\def\CCC{{\mathcal C}}

\def\FFF{{\mathcal F}}

\def\LLL{{\mathcal L}}
\def\MMM{{\mathcal M}}

\def\MMM{{\mathcal M}}

\def\NNN{{\mathcal N}}

\def\SSS{{\mathcal S}}

\def\ZZZ{{\mathcal Z}}

\newtheorem{thm}{Theorem}[section]

\newtheorem{lemma}[thm]{Lemma}
\newtheorem{cor}[thm]{Corollary}

\newcommand{\qed}{\nopagebreak \begin{flushright}
        \rule{2mm}{2.5mm} \end{flushright}}
\newcommand{\qedspecial}[1]{\nopagebreak \begin{flushright}
        \rule{2mm}{2.5mm}{\bf #1} \end{flushright}}
\newcommand{\bdry}{\partial}                     
\newcommand{\id}{\mbox{\rm id}}                  
\newcommand{\cl}{\overline}                      
\newcommand{\Imag}{\mbox{\rm Im}}                    

\newcommand{\Aut}{\mbox{\rm Aut}}                        
\newcommand{\Mod}{\mbox{\rm Mod}}    
\newcommand{\intersect}{\cap}                    
\newcommand{\union}{\cup}                        

\newcommand{\mtwo}[4]                            
{\mbox{$\left(\begin{array}{cc}                  
#1 & #2 \\
#3 & #4 
\end{array}
\right)$}}

\newcommand{\dettwo}[4]                          
{\mbox{$\left|\begin{array}{cc}                  
#1 & #2 \\
#3 & #4 
\end{array}
\right|$}}

\newcommand{\pf}{\noindent {\bf Proof: }}

\newcommand{\be}{\begin{enumerate}}
\newcommand{\eb}{\end{enumerate}}
\newcommand{\bi}{\begin{itemize}}
\newcommand{\ib}{\end{itemize}}
\newcommand{\bl}{\begin{list}}
\newcommand{\lb}{\end{list}}

\newcommand{\gap}{\vspace{5pt}}                 

\newcommand{\genby}[1]{\mbox{$\langle #1 \rangle$}}

\newcommand{\me}{ \noindent \textsc{Kevin M. Pilgrim\\  
                                        Dept. of Mathematics, 
                                        Indiana University, Bloomington, IN\\ 
                                        pilgrim@indiana.edu}}

\newcommand{\psibar}{\overline{\psi}}
\usepackage{amsfonts, color}
\usepackage{graphicx}
\usepackage{amssymb}
\usepackage{epstopdf}
\newcommand{\Tw}{\mbox{\rm Tw}}
\newcommand{\pmcg}{\mbox{\rm Mod}}
\newcommand{\wtP}{\widetilde{P}}
\newcommand{\wtU}{\widetilde{U}}
\newcommand{\wtA}{\widetilde{A}}
\newcommand{\tg}{\tilde{g}}
\newcommand{\tgamma}{\tilde{\gamma}}
\newcommand{\tih}{\tilde{h}}
\newcommand{\ta}{\tilde{a}}

\newcommand{\dom}{\mbox{\rm dom}}

\newcommand{\supp}{\mbox{\rm supp}}

\newcommand{\psihat}{\widehat{\psi}}
\newcommand{\pullback}{\stackrel{f}{\longleftarrow}}

\title{An algebraic formulation of Thurston's characterization of rational functions}
\author{Kevin M. Pilgrim}
\date{\today}                                           

\begin{document}
\maketitle

\begin{abstract}
Following Douady-Hubbard and Bartholdi-Nekrashevych, we give an algebraic formulation of Thurston's characterization of rational functions.   The techniques developed are applied to the analysis of the dynamics on the set of free homotopy classes of simple closed curves induced by a rational function.  The resulting finiteness results yield new information on the global dynamics of the pullback map on Teichm\"uller space used in the proof of the characterization theorem.    
\end{abstract}

\tableofcontents

\newpage

\section{Introduction}
{\em Thurston maps} $f: S^2 \to S^2$ are topological objects, regarded up to homotopy, which arise in the classification and characterization of certain holomorphic dynamical systems on the Riemann sphere \cite{DH1}.  They are two-real-dimensional generalizations of the {\em kneading data} associated to interval maps, introduced by Milnor and Thurston \cite{milnor:thurston}.   

In this work, we show that newly introduced  algebraic invariants can be used to reformulate Thurston's original characterization theorem (Theorem \ref{thm:Thurston_characterization} below), and that this connection can be applied to study a previously mysterious problem, namely, the fate of curves under iterated pullback.  When combined with a recent observation of Selinger \cite{selinger:wp}, this yields detailed results about the global dynamics of Thurston's pullback self-map on Teichm\"uller space.  This pullback map plays a key role  in the proof of Theorem \ref{thm:Thurston_characterization}.  Like its cousin the {\em skinning map}---another self-map of Teichm\"uller space, employed in the hyperbolization of 3-manifolds---it has been the subject of recent detailed investigations; cf. \cite{bekp}, \cite{kent:skinning}, and recent work of D. Dumas.

We briefly review some fundamental concepts from \cite{DH1}.  
\gap

\noindent{\bf Thurston maps.} Let $S^2$ denote the two-sphere equipped with an orientation, let $f: S^2 \to S^2$ be an orientation-preserving branched covering map
of degree $d \geq 2$, and let $P_f = \cup_{n>0}f^{\circ n} (\Omega_f)$
where $\Omega_f$ is the set of branch points of $f$.  We say $f$ is a {\em Thurston map} if $P_f$
is finite; we assume this throughout this work.  Two Thurston
maps $f,g$ are {\em equivalent} if there are
homeomorphisms $h, \widetilde{h}: (S^2, P_f) \to (S^2, P_g)$ such that
$h$ is isotopic to $\widetilde{h}$ relative to $P_f$ and $h \circ f
= g \circ \widetilde{h}$.  The {\em orbifold
$\mathcal{O}_f$ associated to $f$} is the topological orbifold whose 
underlying space $S^2$ and whose weight $\nu(x)$ at $x$ is given by $\nu(x) = \mbox{\rm lcm}\{\deg(f^k, (\tilde{x}) : f^k(\tilde{x})=x\}$; here $\deg( , )$ denotes local degree.  The orbifold $\mathcal{O}_f$ is said to be {\em hyperbolic} if
the Euler characteristic $\chi(\mathcal{O}_f)=2-\sum_{x \in P_f}
(1-1/\nu(x))$ is negative.
\gap

In \cite[Appendix, Example 1]{DH1}, Douady and Hubbard posed the following problem.  Consider $f(z)=z^2+i$; then $P_f = \{\infty, i, i-1, -i\}$.  Let $D$ be a Dehn twist about a simple closed curve in $\rs - P_f$.  Then the branched covering obtained by postcomposing $f$ with $D$ is either equivalent to $z^2+i$, to $z^2-i$, or else is not equivalent to a quadratic polynomial; the problem is to characterize the outcome as a function of $D$, regarded as an element of the pure mapping class group of homeomorphisms of the sphere fixing $P_f$ pointwise.  The analogous problem for the so-called ``rabbit'' polynomial became known as the {\em Twisted Rabbit Problem}.

For some years, these problems remained unsolved and were a humbling reminder of the lack of a suitable arsenal of invariants.  The situation changed with the publication  of \cite{bartholdi:nekrashevych:twisted}, in which the newly developed theory of {\em selfsimilar groups} was brought to bear to give complete solutions to these problems.    In addition, in the presence of suitable expansion, Thurston maps admit cellular Markov partitions \cite{bonk:meyer:subdivisions}.  The result is that the combinatorial theory of Thurston maps is now immensely richer; see e.g. \cite{nekrashevych:combinatorics}, \cite{kelsey:schemes}.  

The concept of a {\em virtual endomorphism} plays a key role in the solution to the Twisted Rabbit Problem.  
\gap

\noindent{\bf Virtual endomorphism.}  If $G$ is a group and $H < G$ is a subgroup of finite index, a homomorphism $\phi: H \to G$ is called a {\em virtual endomorphism} of $G$; this is sometimes written $\phi: G \dasharrow G$ and we write $\dom(\phi)=H$.  

Throughout, we deal with finitely generated groups.  Let $|g|$ denote the word length of $g$ with respect to a symmetric generating set $S$ of $G$.  The {\em contraction coefficient} of $\phi$  
\[ \rho(\phi) = \limsup_{n \to \infty} \left( \limsup_{g \in \dom(\phi^{\circ n}), |g|\to \infty} \frac{|\phi^{\circ n}(g)|}{|g|}\right)^{1/n} \]
is independent of the generating set.   One says $\phi$ is {\em contracting} on $G$ if $\rho(\phi)<1$.  

Suppose now that $\phi: G \dasharrow G$ is a virtual endomorphism.   A subgroup $H < G$ is {\em quasi-invariant} if $\phi|_H: H \dasharrow H$ is again a virtual endomorphism, that is, $[H: H\intersect \dom(\phi)]<\infty$ and $\phi(H \intersect \dom(\phi))\subset H$.  
\gap

\noindent{\bf Curves.}   For $P \subset S^2$ a finite set with at least three points,  denote by $\CCC(S^2, P)$ the set of free homotopy classes of essential, unoriented, simple, closed, nonperipheral (that is, not homotopic into arbitrarily small neighborhoods of elements of $P$) curves in $S^2-P$; we use the term {\em curve} for an element of $\CCC(S^2,P)$.  When $P$ is understood we write simply $\CCC$.  A {\em multicurve} $\Gamma$ is a nonempty set of distinct elements of $\CCC(S^2, P)$ represented by pairwise nonintersecting curves.   The set of nonempty multicurves is denoted $\MMM\CCC(S^2, P)$ or simply by $\MMM\CCC$.  

Now suppose $f$ is a Thurston map, and let $\CCC = \CCC(S^2, P_f)$.  Let $o$ denote the union of the homotopy classes in $S^2-P_f$ of curves which are either inessential or peripheral; we call such curves {\em trivial}.  The {\em pullback relation} $\pullback$ on $\CCC \union \{o\}$ is defined by setting $o \pullback o$ and 
\[ \gamma_1 \pullback \gamma_2\]
if and only if  $\gamma_2$ is homotopic in $S^2-P_f$ to a component of the preimage of $\gamma_1$ under $f$.  Thus $\gamma \pullback o$ if and only if some preimage of $\gamma$ is inessential or peripheral in $S^2-P_f$.   The pullback relation induces a {\em pullback function} $f^{-1}: \MMM\CCC \union \{o\} \to \MMM\CCC \union \{o\}$, where $o$ stands for the empty multicurve, by sending $o \mapsto o$ and $\Gamma \mapsto f^{-1}(\Gamma):=\{ \tilde{\gamma} : \exists \gamma \in \Gamma, \gamma \pullback \tilde{\gamma}\}$; here we set $f^{-1}(\Gamma) = o$ if the set of such $\tilde{\gamma}$ is empty.  A multicurve is {\em invariant} if $f^{-1}(\Gamma) \subset \Gamma$ or $f^{-1}(\Gamma) = o$; it is {\em completely invariant} if $f^{-1}(\Gamma) = \Gamma$.  
We say that the pullback function has a {\em finite global attractor} if there exists a finite subset $\NNN \subset \MMM\CCC \union \{o\}$ such that $f^{-1}(\NNN) \subset \NNN$ and for all $\Gamma \in \MMM\CCC$, there exists a nonnegative integer $n$ such that $(f^{-1})^{\circ n}(\Gamma) \in \NNN$; the smallest such subset, if it exists, consists of periodic cycles and is called {\em the} finite global attractor of the function $f^{-1}$.  We extend this concept to the pullback relation in the obvious way.

Denote by $\Z[\CCC(S^2, P)]$ and $\R[\CCC(S^2, P)]$ the free $\Z$- and $\R$-modules generated by $\CCC(S^2,P)$, so that an element $w$ of e.g. $\Z[\CCC(S^2, P)]$ is given by a formal finite linear combination $w=\sum_i a_i \gamma_i$, $a_i \in \Z$.  The free submodules generated by a multicurve $\Gamma$ will be denoted $\Z^\Gamma$ and $\R^\Gamma$.

The {\em Thurston linear transformation} 
\[  \LLL_f: \R[\CCC] \to \R[\CCC]\]
is defined on basis vectors by 
\[ \LLL_f(\gamma) = \sum_{\gamma \pullback \gamma_i} d_i \gamma_i\]
where 
\[ d_i = \sum_{f^{-1}(\gamma)\supset \delta \simeq \gamma_i}\frac{1}{\deg(f:\delta \to \gamma)};\]
the sum ranges over preimages $\delta$ of $\gamma$ homotopic to $\gamma_i$.  

{\em Thurston's Characterization Theorem} \cite[Theorem 1]{DH1} asserts 

\begin{thm}[Thurston's characterization]
\label{thm:Thurston_characterization}
If $\mathcal{O}_f$ is hyperbolic, then $f$ is equivalent to a rational map $R$ if and only if for every invariant multicurve $\Gamma$, the spectrum of the linear map $\LLL_{f,\Gamma}$ lies strictly inside the unit disk; in this case, $R$ is unique, up to conjugation by M\"obius transformations.   
\end{thm}

Though a relation, not a function, it is natural to ask the usual dynamical questions regarding iteration of the pullback relation $\pullback$: can curves be periodic? how many cycles can exist? can curves wander? It is easy to construct Thurston maps $f$ with invariant subsurfaces $\Sigma \subset S^2$ on which the map $f$ acts like an arbitrary element of the mapping class group of $\Sigma$; such maps are necessarily obstructed.  Hence questions about iteration of $\pullback$ are most naturally posed in the setting when $f$ is a rational Thurston map.  Here, tension arises.  On the one hand, if $\gamma \pullback \tilde{\gamma}$ and $\deg(f:\tilde{\gamma} \to \gamma)>1$, then the unique hyperbolic geodesic in $\rs \setminus P_f$ homotopic to $\gamma$ lifts under $f$ to a longer hyperbolic geodesic in $\rs \setminus f^{-1}(P_f)$ and so $\tilde{\gamma}$ might be more complicated than $\gamma$.  On the other hand, the inclusion $\rs \setminus f^{-1}(P_f) \hookrightarrow \rs \setminus P_f$ contracts hyperbolic metrics.  Hence {\em a priori}, it is unclear which phenomenon---length increase or decrease---has the dominant effect.  
\gap

We now connect the pullback relation on curves with a certain virtual endomorphism.  
\gap

\noindent{\bf Mapping class groups.}  For $P \subset S^2$ the {\em pure mapping class group} $\Mod(S^2, P)$ is the group of homotopy classes of orientation-preserving homeomorphisms $h: S^2 \to S^2$ which fix $P$ pointwise.   There is a distinguished subset $\Tw(S^2, P) \subset \Mod(S^2,P)$ consisting of {\em multitwists}, that is, mapping class elements represented by products of powers of Dehn twists about the elements of a multicurve $\Gamma$.   The set $\Tw(S^2, P)$ is invariant under conjugation, but not under arbitrary group automorphisms; if one restricts to the class of automorphisms preserving a peripheral structure around points in $P_f$, then the set of multitwists becomes characteristic.  There is a natural bijection 
\[ \gamma \leftrightarrow T_\gamma\]
between homotopy classes of unoriented simple nonperipheral curves and {\em left} Dehn twists.  More generally, if $\Gamma=\{\gamma_i\}$ is a multicurve and $w=\sum_i a_i \gamma_i \in \Z^\Gamma$  then the correspondence 
\begin{equation}
\label{eqn:correspondence}
 M_w=\prod_i T_{\gamma_i}^{a_i}\ \mapsto \ \sum_i a_i \gamma_i =w
\end{equation}
defines an injection of sets 
\[ \Tw(S^2, P_f) \hookrightarrow \Z[\CCC(S^2, P)] \subset \R[\CCC(S^2, P)].\]
If $\Gamma$ is a multicurve, we denote by $\Tw(\Gamma)$ the subgroup of $G$ generated by Dehn twists about the elements of $\Gamma$; it is free abelian, of rank $\#\Gamma$.  

Suppose now  $f$ is a Thurston map, and put $G=\Mod(S^2, P_f)$.  Elementary covering space theory implies that there is a finite-index subgroup $H<G$ consisting of mapping classes representable by homeomorphisms that lift under $f$ to homeomorphisms which again represent elements of $G$, i.e. which fix the set $P_f$ pointwise.   We obtain a virtual endomorphism 
\[ \phi_f: G \dashrightarrow G\]
such that for representative homeomorphisms, 
\[ h \circ f = f \circ \phi_f(h).\]
\gap

\noindent{\bf Results.}  The following result connects topology and algebra:  

\begin{thm} 
\label{thm:phi_and_LLL}
If $M_w \in \dom(\phi_f)$, then 
\[ \phi_f(M_w) = M_{\LLL_f(w)}\]
where $\LLL_f$ is the Thurston linear transformation.  In particular, $\Gamma$ is an invariant multicurve if and only if $\Tw(\Gamma)$ is a $\phi_f$-quasi-invariant subgroup, where $\Tw(\Gamma) = \genby{T_\gamma: \gamma \in \Gamma} <  G$.    In this case, under the correspondence (\ref{eqn:correspondence}), 
\[ \LLL_{f,\Gamma} = (\phi_f|_{\Z^\Gamma})\otimes \R.\]
\end{thm}

\noindent{\bf Applications I:  algebraic characterization}  
As an application, we give an algebraic version of Thurston's Characterization Theorem:  

\begin{thm}
\label{thm:rational_iff_contracting} A Thurston map $f$ with hyperbolic orbifold and with $\#P_f \geq 4$ is equivalent to a rational map if and only if for every $\phi_f$-quasi-invariant abelian subgroup $H < \Tw(S^2, P_f)$, the induced virtual endomorphism $\phi_f|_H: H \dasharrow H$ has contraction coefficient $<1$.  
\end{thm}

\noindent{\bf Remarks.}

\be
\item If $\#P_f < 3$ then $f$ is equivalent to $z \mapsto z^k$ for some $k$ with $|k| \geq 2$; if $\#P_f < 4$ then $f$ is always equivalent to a rational  map.  We exclude these cases throughout this work to keep the statements clean.  

\item For fixed degree $\deg(f)$ and cardinality $\#P_f$, there are only finitely many possibilities for the matrices of $\LLL_{f,\Gamma}$.  Hence there exists a constant $C=C(\deg(f), \#P_f)<1$ such that if $f$ is rational,  then $\rho(\phi_f|_H: H \dasharrow H) < C$ for all quasi-invariant abelian twist subgroups.   

\item If $f$ is rational, the contraction on quasi-invariant abelian twist subgroups cannot, in general, be extended to contraction on all of $G$; see \S 9 below.  
\eb
\gap

\noindent{\bf Applications II:  pullback relation on simple closed curves.}

Fix a Thurston map $f$ and let $\TTT$ be the Teichm\"uller space modelled on $(S^2, P_f)$ as in \cite{DH1}.  Associated to $f$ is Thurston's {\em pullback map} $\sigma_f: \TTT \to \TTT$.  Since $\TTT$ is homeomorphic to an open ball and $\sigma_f$ is distance nonincreasing with respect to the Teichm\"uller metric, the dynamics of $\sigma_f: \TTT \to \TTT$ is uninteresting.  However, 
Selinger \cite{selinger:wp} showed that the pullback map $\sigma_f$ extends to the Weil-Petersson completion $\cl{\TTT}$.  The following facts are known; see \cite{wolpert:handbook}.  The completion $\cl{\TTT}$ coincides with the so-called {\em augmented Teichm\"uller space}.  This is a noncompact stratified space whose strata $\TTT_\Gamma$ are in bijective correspondence with multicurves.  Each stratum $\TTT_\Gamma$ is homeomorphic to the product of the Teichm\"uller spaces of the components of the noded surfaces obtained by collapsing exactly those elements of $\Gamma$ to points.  By a theorem of Brock and Margalit \cite{brock:margalit:pants}, $\cl{\TTT}$ is quasi-isometric to the {\em pants complex}.    The definitions immediately imply that 
\[ \sigma_f: \TTT_\Gamma \to \TTT_{f^{-1}(\Gamma)}\]
and so the orbit of a stratum under $\sigma_f$ is encoded by the pullback function $f^{-1}: \MMM\CCC \union \{o\} \to \MMM\CCC\union \{o\}$.  In particular, proper strata invariant under $f^{-1}$ are in bijective correspondence with completely invariant multicurves.   Thus, the extension of $\sigma_f$ to $\cl{\TTT}$ can have interesting dynamics.  

\begin{thm}
\label{thm:contracting implies eventually cycles}
If $\phi_f$ is contracting, then the pullback function on multicurves has a finite global attractor.
\end{thm}

It follows that the pullback relation on curves has a finite global attractor as well.  Note that the converse need not hold; see \S 9.  One can also give analytic conditions on $\sigma_f$ which imply that the pullback function on multicurves has a finite global attractor; this is the subject of ongoing work.  

Using results of Koch and Nekrashevych, we deduce (Corollary \ref{cor:finitely many multicurves for periodic quadratics}) that for critically periodic quadratic polynomials, the pullback function on multicurves has a finite global attractor .  

For general rational maps, we have the following weaker statement, whose proof uses the combination and decomposition theory developed in \cite{kmp:cds}:

\begin{thm}
\label{thm:finitely_many_invariant_multicurves}
Suppose $f$ is a rational map with hyperbolic orbifold.  Then there are only finitely many completely invariant multicurves.  
\end{thm}

The proof of Theorem \ref{thm:contracting implies eventually cycles} suggests a general method to calculate the finite global attractor of the pullback function on multicurves.  We apply this analysis to each of the three quadratic polynomial examples $z \mapsto z^2+c$ with three finite postcritical points studied in \cite{bartholdi:nekrashevych:twisted}.  Among experts in complex dynamics, they are referred to by a rational number mod $1$ known as an {\em external angle}: $f_{1/7}$, also known as the ``rabbit'' polynomial; $f_{1/6}$, the ``dendrite'' polynomial given by the formula $f_{1/6}(z)=z^2+i$; and $f_{1/4}$.   

For the rabbit polynomial $f_{1/7}$, the virtual endomorphism $\phi_f$ is contracting.  We exploit this to prove

\begin{thm}
\label{thm:nwc_for_rabbit}
Let $f=f_{1/7}$.  Under backward iteration, every curve becomes either trivial, or falls into the unique three-cycle.
\end{thm}

In contrast, for the dendrite polynomial $f_{1/4}$, the virtual endomorphism $\phi_f$  is not contracting.  Nevertheless, modified methods yield: 

\begin{thm}
\label{thm:eventually_trivial_dendrite}
Let $f=f_{1/6}$ be the dendrite polynomial.  Under backward iteration, every  curve becomes trivial.
\end{thm}

For the polynomial $f_{1/4}$, the virtual endomorphism $\phi_f$ is again contracting.  R. Lodge, using similar methods as in the proof of Theorem \ref{thm:nwc_for_rabbit}, shows:

\begin{thm}
\label{thm:eventually_trivial_pre1per2}
Let $f=f_{1/4}$.  Under backward iteration, every curve becomes trivial.  
\end{thm}

\gap

\noindent{\bf Organization.}  In \S 2 we discuss in detail mapping class groups, the correspondence $M_w \leftrightarrow w$, and state Thurston's classification for pure mapping classes on the sphere.  We factor the virtual endomorphism $\phi_f$, defined in \S 5,  as a composition of lifting and filling in punctures, discussed respectively in  \S\S 3 and 4.  The proofs of Theorems \ref{thm:phi_and_LLL} and \ref{thm:rational_iff_contracting} are given in \S\S 5 and 6, respectively.  In \S 7 we prove Theorems \ref{thm:contracting implies eventually cycles} and \ref{thm:finitely_many_invariant_multicurves} while  \S\S 8, 9, 10 give the analysis of the maps $f_{1/7}$, $f_{1/6}$, and $f_{1/4}$, respectively.  

\gap

\noindent{\bf Conventions.}  To avoid uninteresting and special cases, we assume throughout that $\mathcal{O}_f$ is hyperbolic and that $\#P \geq 4$ unless otherwise stated. We follow the notational conventions as in  \cite{bartholdi:nekrashevych:twisted}.  In particular, given transformations $S, T$, the notation $ST$ indicates that $S$ is performed first, then $T$, i.e. their action is a right action.  In long expressions we often distinguish between factors in such products by the symbol $\cdot$, so that $S\cdot T = ST$.  The notation $S \circ T$ indicates that $T$ is performed first, then $S$.  If $g_1, g_2$ are elements of a group, conjugation is given as a right action, so that $g_1^{g_2} = g_2^{-1}g_1g_2$.  If $(g_1, \ldots, g_d)\alpha, (h_1, \ldots, h_d), \beta$ are elements of a wreath product $G^d \rtimes S_d$, their product is given by 
\[ (g_1 h_{1^\alpha}, \ldots, g_d h_{d^\alpha})\alpha\circ \beta\]
where $i^\alpha$ is the image of $i$ under the permutation $\alpha$.  
\gap

\noindent{\bf Acknowledgements.}  I thank Volodymyr Nekrashevych and Laurent Bartholdi for useful conversations.

\section{Mapping class groups}

We begin with generalities.  Throughout this section, $P \subset S^2$ is a finite set with at least four points.

Elements of $\pmcg(S^2, P)$ are special.  

\begin{lemma}
\label{lemma:fixed}
Suppose $g \in \pmcg(S^2, P)$ permutes the elements of a multicurve $\Gamma$ up to isotopy.  Then $g$ fixes the elements of $\Gamma$ up to isotopy.  If in addition the elements of $\Gamma$ are oriented, $g$ preserves the orientation of each element.  
\end{lemma}

\pf If $\Gamma$ has only one element $\gamma$ then since $g|P=\id_P$, $g$ preserves each of the complementary components of $\gamma$ up to isotopy and so preserves an orientation on $\gamma$.   We now induct on $\#\Gamma$.     There exists $\gamma \in \Gamma$ which does not separate any pair of elements of $\Gamma$.   Thus there exists a Jordan domain   $D \subset S^2-\Gamma$ bounded by an element $\gamma$ of $\Gamma$ such that $D \intersect P \neq \emptyset$.  Since $g|_P=\id_P$ we must have $g(D)=D$ up to isotopy fixing $P$  and so in particular $g(\gamma)$ is isotopic to $\gamma$ and $g|\gamma$ preserves an orientation on $\gamma$.  Thus $g(\Gamma')=\Gamma'$ where $\Gamma'=\Gamma - \{\gamma\}$.  By induction, the proof is complete.\qed

\noindent{\bf Thurston's classification.}  Thurston's classification of mapping classes  \cite[Theorem 4]{thurston:bulletin:surface_diffeos} is correspondingly simpler.  

\begin{thm}[Thurston classification]
\label{thm:Thurston classification}
A nontrivial element $g \in \pmcg(S^2,  P)$ is either 
\be
\item aperiodic reducible:  of infinite order and permutes (hence fixes) the elements of a nonempty multicurve $\Gamma$,  preserving orientation, or 
\item pseudo-Anosov.   
\eb
\end{thm}

The finite order case cannot occur, since by a classical result of Nielsen \cite{nielsen:realization_for_cyclic}, any periodic mapping class is represented by a M\"obius transformation, and we are assuming that such classes fix each of the $\geq 3$ elements of the set  $P$ of marked points.  
\gap

\noindent{\bf Support.}  Given a homeomorphism $h: (S^2, P) \to (S^2, P)$, its {\em support} $\supp(h)$ is the closure of complement of the set of fixed points of $h$.   The {\em support} of a weighted multicurve $w=\sum_iw_i\gamma_i$ is $\Gamma_w = \union_{w_i\neq 0}\gamma_i$.   
\gap

\noindent{\bf Twists.}  Let $\gamma$ represent an essential simple closed curve and let $A$ be a closed regular neighborhood of $\gamma$ in $S^2-P$.  Then there is an orientation-preserving homeomorphism $\phi: \{1 \leq |z| \leq 2\} \to A$.  A  {\em positive, or left, Dehn twist $T_\gamma$ about $\gamma$} is an element of $\pmcg(S^2, P)$ represented by a homeomorphism $h$ whose support lies in $A$ and which on $A$ is given by $\phi\circ \hat{h}\circ \phi^{-1}$, where $\hat{h}(re^{2\pi i\theta}) = re^{2\pi i(\theta + (r-1))}$.   Thus, if $\alpha$ is the image under $\phi$ of the segment joining $1$ and $2$, then $h(\alpha)$ bends to the left as it is traversed in either direction and winds once around the annulus $A$.   The class of $T_\gamma$ depends only on the isotopy class of $\gamma$.  
  
\begin{lemma}
\label{lemma:compute_p}
Suppose $g=h^p$ where $h$ is a left Dehn twist and $p \in \Z$.  Then $p$ may be computed as follows.  Represent $h$ by a homeomorphism supported on an annulus $A$.  Let $\alpha:([0,1], \{0,1\}) \to (A, \bdry A)$ be a path joining a point $a$ in one boundary component $\gamma$ of $A$ to a point in the other boundary component.  Give $\gamma$ the orientation induced from $A$; then $[\gamma] \mapsto 1$ identifies $\pi_1(A,a)$ with $\Z$.  Then $p=[\alpha * g(\cl{\alpha})]$ where $*$ denotes concatenation of paths, $\cl{\alpha}$ denotes the path $\alpha$ traversed in the opposite direction, and $\alpha$ is traversed first.  
\end{lemma}
 \gap

\noindent{\bf Multitwists.}  A {\em multitwist} is an element of $\pmcg(S^2, P)$ which is a product of powers of Dehn twists about the elements of a multicurve.  We denote by $\Tw(S^2, P)$ the subset of $\pmcg(S^2, P)$ given by multitwists, and by $\Tw^+(S^2, P)$ the subset of twists in which the powers are all strictly positive.  The implied representation as a product is unique:

\begin{lemma}
\label{lemma:multitwists}
Suppose $\{\alpha_1, \ldots, \alpha_m\}$ and $\{\beta_1, \ldots, \beta_n\}$ are multicurves.  Set $a_i=T_{\alpha_i}, b_j = T_{\beta_j}$,  and suppose $p_1, \ldots, p_m, q_1, \ldots, q_n  \in \Z \setminus \{0\}$.  If 
\[g = a_1^{p_1} \ldots a_m ^{p_m} = b_1^{q_1} \ldots b_n^{q_n}\]
then $m=n$ and, after re-indexing if needed,  $a_i = b_i$ and $p_i = q_i$, $i=1, \ldots, n$.  
\end{lemma}

\pf By assumption, $g(\alpha_i)=\alpha_i$ and $g(\beta_j)=\beta_j$ up to isotopy for all $i, j$.  Let $\iota$ denote the geometric intersection number of a pair of simple closed curves, i.e. the infimum of the number of intersection points as the representatives for the classes vary.   If $ \iota(\alpha_i, \beta_j) \neq 0$ for some $i, j$ then  by \cite[Prop. 2.2]{farb:margalit:mcg} we have 
\[  \iota(g^l(\alpha_i), \beta) \to \infty\]
as $l \to \infty$, which is impossible since $g^l(\alpha_i)=\alpha_i$ for each $l$.  Hence the $\alpha_i$'s and the $\beta_j$'s are pairwise disjoint up to homotopy.  The $a_i$'s and the $b_j$'s then commute and together freely generate an abelian subgroup.  In such a group the representation of an element as a product of the given free abelian generators is unique.   
\qed

Thus if $g \in \pmcg(S^2, P)$ is a multitwist, its {\em support}, defined as the multicurve (up to isotopy) about which the nontrivial twists occur, is well-defined.   More generally, the support of a set or group of multitwists is defined as the union of the supports of its elements.  

The subsets $\Tw(S^2,  P)$ and $\Tw^+(S^2,P)$ are  invariant under the action of $\pmcg(S^2, P)$ on itself by conjugation.  
\gap

\noindent{\bf Action on multicurves.}  Since $\pmcg(S^2,P)$ acts on $\CCC(S^2,P)$, there is a representation
\[ \pmcg(S^2,P) \to GL(\R[\CCC(S^2, P)])\]
given in the obvious way by 
\[ g.w = g.\left(\sum_i w_i\gamma_i\right) = \sum_i w_i g(\gamma_i)\]
i.e. by permuting basis elements.  
Note that if $g \in \pmcg(S^2,P)$ then 
\[ g\circ M_w\circ g^{-1} = M_{g.w}.\]
The group $\pmcg(S^2,P)$ also acts on the set of multicurves; there are finitely many orbits.   
\gap

\section{Branched coverings}

Throughout this section, fix a degree $d \geq 2$ and finite subsets $\wtP, P \subset S^2$ where $\#P \geq 4$.  We denote by $\FFF = \FFF(\wtP, P, d)$ the set of homotopy classes of branched coverings $f:(S^2,\wtP) \to (S^2, P)$ of degree $d$ such that $\wtP=f^{-1}(P)$ and $f: S^2\setminus \wtP \to S^2 \setminus P$ is unramified.
\gap

\noindent{\bf Lifts of simple essential curves.}  If $\gamma$ is an essential nonperipheral simple closed curve in $S^2-P$, and if $\{\tgamma_k\}$ denote the  components of $f^{-1}(\gamma)$, then each $\gamma_k$ is again essential and nonperipheral in $S^2-\wtP$, and no two distinct elements $\tgamma_i, \tgamma_j$ are homotopic in $S^2-\wtP$; see \cite[Lemma 1.11]{kmp:cds}.  
\gap

The following lemma is an essential technical ingredient of our analysis.  It will imply that multitwists lift under branched coverings to multitwists.  

\begin{lemma}
\label{lemma:support}
Suppose $f$ represents an element of $f$ and $\{\gamma_1, \ldots, \gamma_n\}$ represents a (possibly empty) multicurve $\Gamma$.  Suppose $A_1, \ldots, A_n$ are pairwise disjoint open regular neighborhoods of $\gamma_1, \ldots, \gamma_n$ in $S^2\setminus P$.  
Suppose $g: S^2 \to S^2$ and $\supp(g) \subset A_1 \union \ldots \union A_n$.

If $g \circ f = f \circ \tilde{g}$ and $\tilde{g}|\wtP = \id_{\wtP}$ then $\supp(\tilde{g}) \subset f^{-1}(A_1 \union \ldots \union A_n)$.   
\end{lemma}

\pf Let $\widetilde{\Gamma} = f^{-1}(\Gamma)$; thus $\widetilde{\Gamma}$ represents a multicurve in $S^2\setminus \wtP$.  By construction $g(\Gamma)=\Gamma$ up to isotopy and so by Lemma \ref{lemma:fixed} up to isotopy $g$ fixes each $\gamma_i \in \Gamma$ and preserves its orientation.  
The same reasoning applied to $\tilde{g}$ implies that up to isotopy $\tilde{g}$ fixes and preserves the orientation of each element of $\widetilde{\Gamma}$.  Now let $U$ be a component of $S^2 \setminus (A_1 \union \ldots \union A_n \union P)$.  The boundary $\bdry U$ consists of punctures and curves isotopic to elements of $\Gamma$, so since $\Gamma$ is a multicurve and $\#P \geq 4$ the Euler characteristic $\chi(U)$ satisfies $\chi(U) \geq 3$.  Upstairs, the boundary $\bdry \wtU$ consists of punctures and curves isotopic to elements of $\widetilde{\Gamma}$.  It follows that if $\widetilde{U}$ is a component of $f^{-1}(U)$ then $\widetilde{g}$ sends each boundary component of $\widetilde{U}$ to itself.  Since $\chi(\widetilde{U}) = \deg(f:\widetilde{U} \to U)\cdot \chi(U) \geq 3$ the Lefschetz fixed-point formula then implies that $\tilde{g}$ has a fixed-point in $\widetilde{U}$.  
The restriction $f: \widetilde{U} \to U$ is an unramified covering and  by assumption $g|U=\id_U$, so the equation $g\circ f = f \circ \tilde{g}$ implies that the restriction $\tilde{g}|_{\widetilde{U}}$ is a covering automorphism of $f: \widetilde{U} \to U$.   Since $\tilde{g}|_{\widetilde{U}}$ has a fixed-point in $\widetilde{U}$, it must be the identity there.

\qed

The groups $\pmcg(S^2,P)$ and $\pmcg(S^2,\wtP)$ act on $\FFF$  on the left and right, respectively, by 
\[ g.f.\tg = g\circ f \circ \tg.\]
Since it is easily verified that these actions are well-defined, we use the notation $\circ$ to denote these actions.  The left action is free.  To see this, suppose $h \circ f \simeq f$ as elements of $\FFF$.  
Then the induced map $h_*$ on $\pi_1(S^2\setminus P)$ is (up to conjugacy) the identity on the finite-index subgroup $f_*\pi_1(S^2 \setminus \wtP)$.  It follows that $h\simeq \id$ in $\pmcg(S^2, P)$.  Lemma \ref{lemma:support} applied to the empty multicurve implies that the right action of $\pmcg(S^2,\wtP)$ on $\FFF$ is also free.  

Let $f \in \FFF$ and $g \in \pmcg(S^2,P)$.  If $g\circ f = f\circ \tg$ for some $\tg \in \pmcg(S^2,\wtP)$ then by freeness of the right action, the element $\tg$ is unique; we denote it $f^*(g)$.  Since the left action is free, $f^*$ is injective.  The set 
\[ \dom(f^*) = \{g | \exists \tg \;\; \mbox{with}\;\; g\circ f=f\circ \tg\}\]
is a subgroup of finite index in $\pmcg(S^2,P)$, and $f^*: \dom(f^*) \to \pmcg(S^2,\wtP)$ is an injective homomorphism.  
\gap

The previous lemma implies 

\begin{lemma}
\label{lemma:lifts_of_twists}
\be
\item 
The homomorphism $f^*$ has the property that
\[ f^*: \dom(f^*) \intersect \Tw(S^2,P) \to \Tw(S^2, \wtP)\]
and
\[ f^*: \dom(f^*) \intersect \Tw^+(S^2,P) \to \Tw^+(S^2,\wtP).\]
In particular, $f^*$ preserves the property of being aperiodic reducible.
\item If $\Gamma = \{\gamma_j\}$ is a multicurve in $\MMM\CCC(S^2, P)$, and $w=\sum_j w_j\gamma_j$, $w_j \in \Z$, then the multitwist $M_w \in \dom(f^*) \iff T_{\gamma_j}^{w_j} \in \dom(f^*) $ for each $j$.
\eb
\end{lemma}

We now refine this observation.  Define a linear transformation 
\[ f^\dagger: \R[\CCC(S^2, P)] \to \R[\CCC(S^2,\wtP)]\]
by 
\[ f^\dagger(\gamma) = \sum_k \frac{1}{d_k}\tgamma_k\]
where $\{\tgamma_k\}$ is the set of components of  $f^{-1}(\gamma)$ and $d_k=\deg(f:\tgamma_k \to \gamma)$ is the (positive) degree.  Then $f^\dagger$, as a linear transformation, depends only on the homotopy class of $f$.  
 
\begin{lemma}
\label{lemma:exponents}
Suppose $g=h^p$ where $h$ is a left Dehn twist about $\gamma$.
\be
\item If $f^\dagger(\gamma)=\sum_k\frac{1}{d_k}\tgamma_k$ and $d_k | p$ for each $k$, 
then $g\in\dom(f^*)$ and $f^*(g) = \prod_k \ta_k^{p/d_k}$
where $\tilde{a}_k$ is the left Dehn twist about $\tgamma_k$.
\item Conversely, if $g\in\dom(f^*)$ and $f^*(g) = \prod_k \ta_k^{q_k}$
where $\tilde{a}_k$ is the left Dehn twist about $\tgamma_k$, then $f^\dagger(\gamma)=\sum_k\frac{1}{d_k}\tgamma_k$ and $p=q_kd_k$ for each $k$.
\eb
\end{lemma} 

\pf Represent $h$ by an element supported on an annulus $A$ and let $\{\wtA_k\}$ denote the components of $f^{-1}(A)$.  

\noindent{\bf 1.}  Note that $d_k = \deg(f: \wtA_k \to A)$.  By Lemma \ref{lemma:support}, the conclusion $g \in \dom(f^*)$ will follow from the existence of an extension of the identity map on $S^2-\union_k\wtA_k$ to a lift $\tilde{g}$ of $g$.  For this to hold, in turn it is enough to check that $h^p|_{A}$ lifts under $f|_{\wtA_k}:\wtA_k \to A$ to a map $\tg_k: \wtA_k \to \wtA_k$ which is the identity on $\bdry\wtA_k$.  The hypothesis that $d_k | p$ , Lemma \ref{lemma:compute_p}, and standard covering space arguments yield the conclusion.  
\gap

\noindent{\bf 2.}  If $g$ lifts under $f$ to a map $\tg$ representing an element of $\pmcg(S^2,\wtP)$ then Lemma \ref{lemma:support} implies that the support of $\tg$ is contained in $\union_k\wtA_k$.  Focusing on a single component $f|_{\wtA_k}:\wtA_k \to A$, Lemma \ref{lemma:compute_p} again implies that $p=q_kd_k$ where $d_k = \deg(f|_{\wtA_k})$.  
\qed

The conclusion of Lemma \ref{lemma:exponents} may be phrased alternatively as follows.
\begin{cor}
\label{cor:exponents} If $\gamma \in \CCC(S^2, P)$ and $p \in \Z$ then 
\[ T_{\gamma}^p \in \dom(f^*) \iff f^\dagger(p\gamma) \in \Z[\CCC(S^2,\wtP)].\]
If $T_{\gamma}^p \in \dom(f^*)$, then 
\[ f^*(T^p_\gamma) = \prod_kT_{\tgamma_k}^{q_k} \iff f^{\dagger}(p\gamma) = \sum_k q_k\tgamma_k.\]
\end{cor}

\section{Forgetful maps}  

In this section, we assume the setup of the previous section.  Now, however, we make the additional assumption that $\wtP \supset P$.  
\gap

The forgetful map $\pi: (S^2,\wtP) \to (S^2, P)$ induces a surjective homomorphism 
\[ \pi_*:  \pmcg(S^2,\wtP) \to \pmcg(S^2,P).\]
It also induces a surjection 
\[ \pi: \CCC(S^2,\wtP) \to \CCC(S^2,P) \union \{o\};\]
those $\tgamma \in \CCC(S^2,\wtP)$ which are inessential or peripheral in $S^2-P$ are sent to $o$.  This in turn induces  a surjective linear map 
\[ \pi_\dagger:  \R[\CCC(S^2,\wtP)] \to \R[\CCC(S^2, P)]\]
defined on basis elements by $\pi_\dagger(\tgamma)=\pi(\tgamma)$ if $\pi(\tgamma)\neq o$ and $\pi_\dagger(\tgamma)=0$ otherwise.  

\begin{lemma}
\label{lemma:forgetful}
The forgetful homomorphism $\pi_*$ sends multitwists to multitwists and, more generally, reducible elements to reducible elements.    
In particular, 
\[ \pi_*(T_{\tgamma}) = T_{\pi(\tgamma)}.\]
Consequently 
\[ \pi_*: \Tw(S^2,\wtP) \to \Tw(S^2,P)\]
and
\[ \pi_*: \Tw^+(S^2,\wtP) \to \Tw^+(S^2,P)\union \{\id \}. \]
These restrictions are surjective.  
\end{lemma}

\begin{cor}
\label{cor:forgetful}
Suppose $\widetilde{\Gamma} = \{\tgamma_k\}$ is a multicurve in $(S^2,\wtP)$ and $\{l_k\}$ are positive integers.  Let $\Gamma = \{\gamma_i\}$ be the multicurve comprised of the nontrivial images of elements of $\widetilde{\Gamma}$ under $\pi$.  

Then 
\[ \pi_\dagger\left(\sum_k l_k\tgamma_k\right) = \sum_i \left(\sum_{\pi(\tgamma_k)=\gamma_i} l_k\right) \gamma_i\]
if and only if
\[ \pi_*\left(\prod_k T_{\tgamma_k}^{l_k}\right) = \prod_i T_{\gamma_i}^{m_i}.\]
where $m_i=\displaystyle{\sum_{\pi(\tgamma_k)=\gamma_i}}l_k$.  
\end{cor}

\pf If $\pi_\dagger(\tgamma_k)=\gamma_i$ then $\pi_*(T_{\tgamma})=\gamma_i$ by Lemma \ref{lemma:forgetful}.  Necessity follows since $\pi_*$ is a homomorphism and $\pi_\dagger$ is linear.  

To prove sufficiency, suppose that 
\[\pi_\dagger\left(\sum_k l_k\tgamma_k\right) = \sum_j b_j \beta_j\]
where the $\beta_j$ are distinct.  
By the definition of $\pi_\dagger$, since the $\tgamma_k$'s are disjoint, the $\beta_j$'s are disjoint.  
Since we have already proved necessity, we have then 
\[ \pi_*\left(\prod_k T_{\tgamma_k}^{l_k}\right) = \prod_i T_{\gamma_i}^{m_i} = \prod_j T_{\beta_j}^{b_j}.\]
By the uniqueness Lemma \ref{lemma:multitwists}, we have after permuting indices that $\beta_j = \gamma_i$ and $b_j = m_i$.   The result follows since the $\beta_j$ are distinct, hence linearly independent, and $\pi_\dagger$ is linear.   
\qed

\section{The virtual endomorphism}

In this section, we assume $f: S^2 \to S^2$ is a Thurston map; we set $P=P_f$ and $\wtP=f^{-1}(P)$.  
As usual we assume $\#P \geq 4$.  Here, we define precisely the virtual endomorphism $\phi_f: \pmcg(S^2, P) \to \pmcg(S^2, P)$, show that it depends only on the homotopy class of $f$ relative to $P$, and prove Theorem \ref{thm:phi_and_LLL}.  
\gap

The following lemma uses the fact that $f(P) \subset P$ and so is a fact about dynamics.  
\begin{lemma}
\label{lemma:unique_lift}
Suppose $f, g$ Thurston maps with common postcritical set $P$, and suppose $f$ and $g$ are homotopic through Thurston maps agreeing on $P$.  
Then there exists a unique homeomorphism $\tih: S^2 \to S^2$ with the following properties:
\be
\item $f=g\circ \tih$, and 
\item $\tih$ is isotopic to the identity through homeomorphisms fixing $P$.  
\eb
\end{lemma}

\pf  {\bf Uniqueness.}  If $\tilde{k}$ is another such homeomorphism then $\tilde{k}^{-1}\circ h$ restricted to $S^2\setminus \wtP$ is an automorphism of the covering space $f:S^2\setminus \wtP \to S^2 \setminus P$ which fixes at least three punctures; the argument given in Lemma \ref{lemma:support} shows that it must be the identity.  
\gap

{\bf Existence.}  Let $I=[0,1]$ denote the unit interval.  Let $F_0=f$ and $F_1=g$ be joined by a homotopy $F:I\times S^2 \to S^2$ such that $F_t|_P=\id_P$ for all $t$.  Let $M=I\times (S^2-P)$ and $\widetilde{M}=\{(t,x) : x \not\in F_t^{-1}(P)\}$.  Then the map $\widetilde{M} \to M$ given by $(t,x) \mapsto (t,F_t(x))$ is a covering.  Fix a basepoint $b \in S^2-P$ and $\tilde{b} \in F_0^{-1}(b)$.  The path $t \mapsto (t,b)$ may be lifted to a path $t \mapsto (t, \tilde{b}_t)$ where $\tilde{b}_t \in F_t^{-1}(b)$.  The images of the fundamental groups
\[ (F_t)_*:\pi_1(S^2-F_t^{-1}(P), \tilde{b}_t) \to \pi_1(S^2,b)\]
are constant and so the identity map $S^2-P \to S^2-P$ lifts under $F_0=f$ and $F_t$ to a map $\tih_t$ such that $f=F_t\circ \tih_t$.  The resulting family of maps $\tih_t$ is continuous in $t$.  By construction $\tih_t|_P=\id_P$.   Taking $\tih=\tih_1$ we see that $f=g\circ \tih$ with $\tih$ isotopic to the identity relative to $P$.  
\qed

\noindent{\bf Associated virtual endomorphism.} 
Now suppose $f: (S^2, \wtP_f) \to (S^2, P)$ and $g: (S^2, \wtP_g) \to (S^2, P)$ are two Thurston maps with common postcritical set $P$ which are homotopic through Thurston maps agreeing on $P$.  Let $\tih: (S^2, \wtP_f) \to (S^2, \wtP_g)$ be the canonical homeomorphism given by Lemma \ref{lemma:unique_lift} and let $f^*, g^*$ denote the virtual homomorphisms of \S 3.  If $T$ represents an element of $\dom(f^*)$ which lifts under $f$ to a map $\tilde{T}$ fixing $\wtP_f$ then $\tih \circ \tilde{T}\circ \tih^{-1}$ is a lift of $T$ under $g$ which fixes $\wtP_g$.  Let $\tih_*:  \pmcg(S^2,\wtP_f) \to \pmcg(S^2, \wtP_g)$ denote the isomorphism induced by $\tih$.  Then $h_* \circ f^* = g^*$.  
From \S 4, let $\pi^f$, $\pi^g$ denote the forgetful maps from $(S^2, \wtP_f), (S^2, \wtP_g)$ to $(S^2,P)$, respectively.   Then $\pi_*^g \circ h_* = \pi_*^f$.  In summary, the compositions satisfy 
\[ \pi_*^f \circ f^* = \pi_*^g\circ g^*.\]
Thus, given a Thurston map $f$,  the {\em associated virtual endomorphism of $\pmcg(S^2,P)$}  given by the composition $\phi_f = \pi_*^f \circ  f^*$ depends only on the homotopy class of $f$ relative to $P$.  
\gap

Lemmas \ref{lemma:forgetful} and \ref{lemma:lifts_of_twists} imply immediately that 
\[ \phi_f: \dom(\phi_f) \intersect \Tw(S^2,P)  \to \Tw(S^2,P)\]
and
\[ \phi_f: \dom(\phi_f) \intersect \Tw^+(S^2,P) \to \Tw^+(S^2,P) \union \{\id\}.\]
\gap

Also,  the Thurston linear transformation
\[ \LLL_f: \R[\CCC(S^2, P)] \to \R[\CCC(S^2, P)]\]
factors as the composition 
\[ \LLL_f = \pi_\dagger \circ f^\dagger .\]
\gap

\noindent{\bf Proof of Theorem \ref{thm:phi_and_LLL}.}
Suppose $\Gamma = \{\gamma_1, \ldots, \gamma_n\}$ is a multicurve in $\MMM\CCC(S^2, P)$ and $w=\sum_j w_j \gamma_j$, $w_j \in \Z$, is a weighted multicurve.  By Lemma \ref{lemma:lifts_of_twists}, $M_w \in \dom(\phi_f) \iff T_j^{w_j} \in \dom(\phi_f)$ for each $j$.  Fix $j$ and set $\gamma = \gamma_j$ and $p=w_j$.  

By the definitions, Corollary \ref{cor:forgetful}, and Corollary \ref{cor:exponents}, we have   
\[T_\gamma^p \in \dom(\phi_f) \iff \LLL_f(p\gamma) \in \Z[\CCC(S^2, P)]\]  
and in this case, 
\[ \phi_f(T_\gamma^p) = \prod_i T_{\gamma_i}^{l_i} \iff \LLL_f(p\gamma) = \sum_i l_i \gamma_i.\]
Since $\phi_f$ is a homomorphism and $\LLL_f$ is linear, it follows that 
\[ M_w \in \dom(\phi_f) \implies \phi_f(M_w) = M_{\LLL_f (w)}.\]

Now suppose $\Gamma$ is an invariant multicurve, and consider $\LLL_{f,\Gamma}: \R^\Gamma \to \R^\Gamma$.   The definition of $\LLL_f$ implies that there exists a positive integer $L$ such that $\LLL_{f,\Gamma}(L\gamma) \in \Z^\Gamma$ for all $\gamma \in \Gamma$.   Let $H=\genby{T_\gamma : \gamma \in \Gamma} < \Tw(S^2, P_f)$ and identify $H$ with $\Z^n$ by means of the generators so that $\LLL_{f,\Gamma} = (\phi_f|H)\otimes \R$.  Since 
\[ (L\Z)^n < \dom(\phi_f|_H) < \Z^n,\]
it follows that the restriction $\phi_f: H \dasharrow H$ is a virtual endomorphism, and so $H$ is $\phi_f$-quasi-invariant.  

The sufficiency follows from similar arguments.  
\qed

\noindent{\bf Remark:}  There virtual endomorphism defined above is the restriction of another, perhaps more  natural virtual endomorphism.  Namely, one may consider the subgroup $H'_f < \Mod(S^2, P_f)$ consisting of those mapping classes represented by homeomorphisms $h'$ for which (i) there exists a lift $\tilde{h}'$ of $h$ under $f$, and (ii) $\tilde{h'}$ fixes the points of $P_f$.  That the correspondence $h \mapsto h'$ descends to a well-defined homomorphism $H' \to \Mod(S^2, P_f)$,  and that it sends twists to twists, is not immediately obvious; the issue is possibility of covering automorphisms.   In many low-complexity cases, such as quadratic polynomials, the two definitions coincide.

\section{Algebraic Thurston's characterization}

In this section, we prove Theorem \ref{thm:rational_iff_contracting}.  We exploit the fact that the linear map $\LLL_f$ is nonnegative.  

\begin{lemma}
\label{lemma:rho_is_lambda}
Let $N \geq 1$, let $\phi: \Z^N \to \Z^N$ be a virtual endomorphism, and let $A=\phi \otimes \R$ be the associated $\R$-linear map.  Suppose $A$ is nonnegative.  Then the contraction coefficient $\rho(\phi)$ is equal to the Perron-Frobenius leading eigenvalue  $\lambda(A)$.
\end{lemma}

\pf For $v=(v_1, v_2, \ldots, v_N) \in \R^N$ denote by $|v|=\sum_i |v_i|$ the $L^1$-norm on $v$, and for a real nonnegative $N$-by-$N$ matrix $A$ let $||A||$ denote the corresponding operator norm.  

By the Perron-Frobenius theorem and the spectral radius formula \cite[Thm. 14.16]{dym:linear_algebra}, respectively, 
\[ \lambda(A) = r_\sigma(A) = \lim_{n \to \infty}||A^n||^{1/n}=\limsup_{n\to \infty}||A^n||^{1/n}\]
where $r_\sigma(A)$ is the spectral radius of $A$.  
By the definition of the operator norm, this in turn equals 
\[ \limsup_{n\to \infty}\left( \limsup_{0 \neq v \in \R^N} \frac{|A^n(v)|}{|v|}\right)^{1/n}\]
which upon approximating by rationals equals  
\[\limsup_{n\to \infty}\left( \limsup_{0 \neq v \in \Q^N} \frac{|A^n(v)|}{|v|}\right)^{1/n}\]
which in turn, upon clearing denominators, equals 
\[ \limsup_{n\to \infty}\left( \limsup_{0 \neq v \in \Z^N} \frac{|A^n(v)|}{|v|}\right)^{1/n}.\]
Since $\dom(\phi^{\circ n})$ has finite index, there is some integer $L_n>0$ such that $L_nv \in \dom(\phi^{\circ n})$ for all $v \in \Z^N$.  By scaling by $L_n$, the above quantity becomes 
\[\limsup_{n\to \infty}\left( \limsup_{v \in \dom(\phi^{\circ n}), |v| \to \infty} \frac{|\phi^{\circ n}(v)|}{|v|}\right)^{1/n}=\rho(\phi)\]
as required.
\qed

\noindent{\bf Proof of Theorem \ref{thm:rational_iff_contracting}.}  By Theorem \ref{thm:Thurston_characterization}, $f$ is equivalent to a rational map if and only if every $f$-invariant multicurve $\Gamma$ satisfies $\lambda(\LLL_{f, \Gamma}) < 1$.  So suppose $\Gamma$ is an invariant multicurve.  By Theorem \ref{thm:phi_and_LLL}, $H=\Tw(\Gamma)=\genby{T_\gamma : \gamma \in \Gamma}$ is a $\phi_f$-quasi-invariant subgroup, and $\LLL_{f,\Gamma}=(\phi_f|_H)\otimes \R$.  But Lemma \ref{lemma:rho_is_lambda} implies 
$\lambda(\LLL_{f, \Gamma})=\rho(\phi_f|_H)$, so $\lambda(\LLL_{f, \Gamma})<1$ if and only if $\rho(\phi_f|_H)<1$.
\qed

\section{Pullback function on multicurves}

In this section, we prove Theorems \ref{thm:contracting implies eventually cycles} and \ref{thm:finitely_many_invariant_multicurves}.

Suppose $f$ is a Thurston map with postcritical set $P$ and denote $G=\pmcg(S^2, P)$.  Let $\phi_f: G \dasharrow G$ be the virtual endomorphism on $G$ determined by $f$, set $H=\dom(\phi_f)$, and let $\{g_i\}_{i \in I}$ be a choice of coset transversal to $H$ in $G$.

Given a multicurve $\Gamma$, the group $\Tw(\Gamma)$ is a free abelian group; for $N \geq 1$ an integer we denote by $N\cdot \Tw(\Gamma) = \{ g^N : g \in \Tw(\Gamma)\}$.  For each $\Gamma \in \MMM\CCC$, there is a smallest positive integer $N=N(\Gamma) \geq 1$ such that $N \cdot \Tw(\Gamma) < H$.   

Recall that there is a left action of $G$ on $\MMM\CCC \union \{o\}$ by $g.\Gamma = g(\Gamma)$ and $g.o = o$.   There are finitely many $G$-orbits, and since $H$ has finite index in $G$, there are finitely many $H$-orbits as well.  Let $\{\Gamma_j\}_{j \in J}$ be an orbit transversal to the action of $H$ on $\MMM\CCC \union \{o\}$.   
Note that since $H$ consists of classes representable by homeomorphisms that lift under $f$, we have $N(\Gamma)=N(h.\Gamma)$ for all $h \in H$.  

Theorem \ref{thm:phi_and_LLL} implies that 
\[ f^{-1}(\Gamma)=\widetilde{\Gamma} \iff \supp \left[\ \phi_f(N(\Gamma)\cdot \Tw(\Gamma))\ \right]  = \widetilde{\Gamma}.\]

For $(i,j) \in I \times J$ let $\Gamma_{(i,j)} = g_i. \Gamma_j$ and let $N_{(i,j)}=N(\Gamma_{(i,j)})$. 
Given $(i,j)$ there are a unique $\nu(i,j) \in J$ and an $h_{(i,j)} \in H$ for which $f^{-1}(\Gamma_{(i,j)}) = h_{(i,j)}.\Gamma_{\nu(i,j)}$.  Combining this observation with the previous paragraphs, we have 
\[ \supp\left[ \phi_f\left(\ N(\Gamma_{(i,j)}) \cdot \Tw(\Gamma_{(i,j)}) \ \right) \right] = \supp \left[ \Tw(\Gamma_{\nu(i,j)})^{h_{(i,j)}} \right] .\]

 Since $J$ is a transversal to the action of $G$ on $\MMM\CCC\union\{o\}$, there is a surjective map
\[ \pi: J \times G \to \MMM\CCC\union\{o\}\]
defined by 
\[ \pi(\Gamma_j, g) = g.\Gamma_j.\]

Define a function 
\[ \mathbf{f}: J \times G \to J \times G\]
as follows: given $(j,g)$, there are unique $i \in I$ and $h \in H$ with $g=g_i\cdot h$.  We set 
\[ \mathbf{f}(j,g) = (\nu(i,j), h_{(i,j)}\cdot \phi_f(h))\]
where $h_{(i,j)} \in H$ and $\nu(i,j) \in J$ are defined as above.  

\begin{thm}
\label{thm:covers_pullback}
We have 
\[ 
\begin{array}{ccc}
J \times G& \stackrel{\mathbf{f}}{\longrightarrow} & J \times G
\\
\pi \downarrow & \; & \downarrow \pi \\
\MMM\CCC \union \{o\} & \stackrel{f^{-1}}{\longrightarrow} &\MMM\CCC\union\{o\}\\
\end{array}.
\] 
That is, $\Gamma \pullback \tilde{\Gamma}$ if and only if for any $(j,g)$ with $\pi(j,g) = \Gamma$, the multicurve $\tilde{\Gamma}$ is given by  $g'(\Gamma_{j'})$, where $\mathbf{f}(j,g)=(j', g')$.  
\end{thm}

\pf Suppose $\Gamma = \pi(j,g)$.  Let $g_i, h$ be the unique elements satisfying $g=g_i\cdot h$.  Recalling our notational convention, this means that when represented by homeomorphisms, $g=h\circ g_i$, so $\Gamma = h(g_i(\Gamma_j))$.  Thus $\Tw(\Gamma) = \Tw(\Gamma_{(i,j)})^h$.  Then $f^{-1}(\Gamma)= \widetilde{\Gamma}$ if and only if  
\renewcommand{\arraystretch}{2}
\[
\begin{array}{cl}
 &  \tilde{\Gamma} = \supp \left[ \phi_f\left(\ N(\Gamma)\cdot \Tw(\Gamma)\ \right) \right]\\
 \iff & \tilde{\Gamma} =\supp \left[ \phi_f\left(\ N(\Gamma_{(i,j)}) \cdot \Tw(\Gamma_{(i,j)})^h\ \right) \right] \\
 \iff & \tilde{\Gamma} = \supp \left[ \Tw(\Gamma_{\nu(i,j)})^{h_{(i,j)} \cdot \phi_f(h)} ) ) \right] \\ 
 \iff & \tilde{\Gamma} = \pi\left(\ \nu(i,j), h_{(i,j)}\cdot \phi_f(h)\ \right) \\
 \iff & \tilde{\Gamma} = \pi(\mathbf{f}(j, g)).
\end{array}
\]
\qed
Thus any orbit of a multicurve $\Gamma_0 \in \pi(j_0, g_0)$ under iteration of the pullback relation lies in the image under $\pi$ of the orbit of $(j_0, g_0)$ under iteration of the function $\mathbf{f}$.
\gap

\noindent{\bf Remark:}  One could write $\mathbf{f}:J \times G \to G \times J$ and regard $\mathbf{f}$ as defining an automaton.  It seems difficult to apply the connections between automata and selfsimilar group actions as developed in \cite{nekrashevych:book:selfsimilar} to $\mathbf{f}$, however, since the corresponding transformations of $J$ are not in general invertible.
\gap

\noindent{\bf Proof of Theorem \ref{thm:contracting implies eventually cycles}.} 
In the course of the proof, we will make free use of facts from the theory of selfsimilar groups; see \cite{nekrashevych:book:selfsimilar}.
\gap

\pf Set $D=[G:H]$.  Let $\Phi: G \to G^D \rtimes \SSS_D$ be any associated wreath recursion defined by \cite[\S 2.5.5, Equation (2.5)]{nekrashevych:book:selfsimilar}.  The image of an element under $\Phi$ is denoted 
\[ \Phi(g) = \langle g|_1, \ldots, g|_d\rangle \tau(g).\]
By construction, $\phi(h)=h|_1$ for all $h \in H$.  Repeated application of the homomorphism $\Phi$ defines a selfsimilar action of the group $G$ on the set of all finite words $v$ in the alphabet $\{1, \ldots, D\}$.  We denote by $|v|$ the word length of such a word.  

The elements $g|_i$, $i=1, \ldots, D$ are called {\em restrictions}, and repeated application of the wreath recursion gives, to any word $v \in \{1, \ldots, D\}^n$, a corresponding element $g|_v$.  The assigment $(g,v) \mapsto g|_v$ satisfies $g|_{vw} = (g|_v)|_w$ and $(g_1g_2)|_v = (g_1|_{g_2(v)})g_2|_v$.  
By {\em op. cit.}, Lemma 2.11.12, since $\phi$ is contracting, so is $\Phi$.   Fix a generating set for $G$ and let $|| \cdot ||$ denote the associated word length function on $G$.  Contraction of $\Phi$ is equivalent to the following.  There exist natural numbers $n_0, C_0$ such that for all $g \in G$ and all $v$ with $|v| \geq n_0$, the restrictions satisfy $||g|_v|| \leq \frac{1}{2} ||g|| + C_0$.

Suppose $(j_0, g_0) \in J \times G$ is arbitrary, and suppose $(j_n, g_n), n=0, 1, 2, \ldots$ is its orbit under $\mathbf{f}$.  Then for each $n \in \N$, there exist $i_n \in I$, $h_n \in H$ , and $k_n \in \{h_{(i,j)} : (i,j) \in I \times J\}$ such that 
\[ g_n = g_{i_n}\cdot h_n \;\; \mbox{ and }\;\; g_{n+1} = k_{n+1}\cdot \phi_f(h_n)\]
i.e.
\[ h_{n+1} = g_{i_{n+1}}^{-1}\cdot k_{n+1} \cdot \phi_f(h_n).\]
Letting $P=\{g_{i'}^{-1}\cdot  h_{(i,j)} : i, i' \in I, j \in J\}$ and setting $p_n = g_{i_{n+1}}^{-1}\cdot k_{n+1}$, we see that the sequence $\{h_n\}$ satisfies 
\[ h_{n+1} = p_n \cdot \phi_f(h_n), \;\;\;\; n=0, 1, 2, \ldots,\]
where each $p_n$ belongs to the finite set $P$.  Set $P_0 = P$ and for $n \geq 1$ set 
\[ P_n = \{ p|_v : p \in P, |v|=n\}.\]
An easy induction argument and the above properties of restrictions imply that for each $n \in \N$, 
\[ h_n \in P_0 \cdot P_1 \cdot \ldots P_n \cdot h_0|_{\underbrace{11\ldots 1}_{n}}.\]
Thus if $n_0, C_0$ are as above, upon setting 
\[ B_{n_0} = \max\{ ||g|| : g \in P_0 \cdot P_1 \cdot \ldots \cdot P_{n_0}\}\]
we see that 
\[ ||h_{n_0}|| \leq B_0 + \frac{1}{2} ||h_0|| + C_0.\]
Letting $C_1 = \max\{ ||g_i|| : i \in I\}$ we conclude that 
\[ ||g_{n_0}|| \leq C_1 + B_0 + \frac{1}{2} ||h_0|| + C_0 = \frac{1}{2}||g_0|| + C_2.\]
Iterating this bound, it follows that 
\[ ||g_n|| \leq 2C_2 + 1 \]
for all $n$ sufficiently large.  Since $\{ g : ||g|| \leq 2C_2+1\}$ is finite, the proof is complete.
\qed

\begin{cor}
\label{cor:finitely many multicurves for periodic quadratics}
Suppose $f(z)=z^2+c$ and the origin is periodic.  Then the pullback function on multicurves (and, hence, the pullback relation on curves) has a finite global attractor.
\end{cor}

\pf Let $\MMM=\MMM(\IP^1, P_f)$ now be the moduli space of injections $P_f \hookrightarrow \IP^1$ modulo conformal automorphisms of $\IP^1$, so that $\iota \sim \iota'$ if there exists $A \in \Aut(\IP^1)$ such that $\iota' = A \circ \iota$.  Then $\MMM$ is naturally a hyperplane complement in $\C\IP^{\#P_f-3}$ with a distinguished basepoint $\circledast$ corresponding to the identity map of $P_f$.   The mapping class group $G$ is then naturally identified with the fundamental group $\pi_1(\MMM, \circledast)$.  The topology of $\MMM$ is carried by a compact subset $K$, i.e. by the complement of a tubular neighborhood of the omitted hyperplanes.

Under the hypothesis of the corollary, the following facts are known; see \cite{koch:thesis}. There is an associated slightly smaller hyperplane complement $\MMM' \subset \MMM$ and a holomorphic map $\omega_f:  \MMM' \to \MMM$ fixing $\circledast$ (the map $\omega_f$ actually extends to a holomorphic endomorphism of $\C\IP^{\#P_f-3}$, though we do not need this fact).   By \cite[Prop. 3.2]{bekp}, $\omega_f$ is a covering map.  The virtual endomorphism $\phi_f$ on the mapping class group $G$ coincides with the virtual endomorphism on $\pi_1(\MMM, \circledast)$ induced by the covering $\omega_f$ \cite[\S 5]{bartholdi:nekrashevych:twisted}.  The moduli space $\MMM$ is Carath\'eodory hyperbolic, and so one can choose arbitrarily large compact sets $K \subset \MMM$ for which $K \hookrightarrow \MMM$ is surjective on $\pi_1$ and so that $\omega_f^{-1}(K) \subset K$.  It follows that path-lifting of loops under $\omega_f$ uniformly contracts the lengths of loops in $K$, and hence by \cite[Thm. 5.5.3]{nekrashevych:book:selfsimilar} that $\phi_f$ is contracting.   The conclusion then follows by Theorem \ref{thm:contracting implies eventually cycles}.
\qed

\noindent{\bf Remark:} In certain cases when the dynamical map $f$ has higher complexity, one has not a holomorphic map $\omega_f: \MMM' \to \MMM$ on moduli space, but a holomorphic correspondence, i.e. a pair of functions $Y: \MMM' \to \MMM$ and $X: \MMM' \to \MMM$ where $Y$ is a finite covering and $X$ is holomorphic.  The virtual endomorphism $X_* \circ Y_*^{-1}$ on the fundamental group of moduli space again coincides with $\phi_f$.  One can give analytic conditions on this correspondence, similar in spirit to those given in Corollary \ref{cor:finitely many multicurves for periodic quadratics}, which again imply that the pullback relation on curves has a finite global attractor.  However, Lodge (personal communication) has found examples of rational maps for which the pullback relation on curves has a finite global attractor, but for which the source of this finiteness is not a consequence of known general algebraic or analytic properties.
\gap

We now turn to the proof of Theorem \ref{thm:finitely_many_invariant_multicurves}.  
First, note that the hypotheses are clearly necessary, albeit for stupid reasons: Lemma \ref{lemma:3flavors} below shows that for the rabbit polynomial, there are infinitely many distinct curves $\gamma$ which pull back to the trivial curve; thus for each, the multicurve $\Gamma:=\{\gamma\}$  is by definition invariant, but not completely invariant.  And any integral Latt\`es example (see \cite{DH1}) has the property that every multicurve is completely invariant; such maps have Euclidean orbifold.  

The spirit of the proof is the Change of Coordinates Principle \cite[\S 1.3]{farb:margalit:mcg}.  Suppose instead we want to prove that on a closed surface of genus $g$, up to homeomorphism there are only finitely many multicurves.  Cut the surface along a multicurve; one obtains finitely many pieces.  By the classification of surfaces, there are only finitely many possibilities for the pieces.  There are only finitely many ways to glue the pieces together.  Thus, given two multicurves $\Gamma_1, \Gamma_2$ which yield the same pieces and gluing data, decomposing along one and then regluing using the data given by the other yields a homeomorphism sending $\Gamma_1 \to \Gamma_2$. 
\gap 

\pf (Theorem \ref{thm:finitely_many_invariant_multicurves}).   Identify $S^2$ with $\rs$.  Suppose $f$ is a rational map with hyperbolic orbifold; set $Q=f^{-1}(P_f)$.  Consider the pullback relation on homotopy classes of essential, unoriented, simple closed curves in $S^2-Q$; it suffices to prove that a map combinatorially equivalent to $f$ has the property that there are only finitely many completely invariant multicurves under this new relation.  

Suppose $\Gamma$ is such a multicurve.  We decompose $f$ along $\Gamma$ as in \cite{kmp:cds}.  Thicken the elements of $\Gamma$ to a family of annuli $\AAA_0$.  Let $\SSS_0$ be the collection of spheres obtained by cutting along elements of $\AAA_0$ and then adding disks (each with a distinguished marked point in its interior) along each boundary component.  Let $\AAA_1$ be the collection of annuli which are preimages of annuli in $\AAA_0$ and which, up to homotopy, are essential subannuli of $\AAA_0$; by adjusting $f$ if needed within its combinatorial class, we may assume $\AAA_1 \subset \AAA_0$ and $\bdry \AAA_1 \supset \bdry \AAA_0$.  Let $\SSS_1$ be the collection of spheres obtained by cutting along elements of $\AAA_1$ and then adding disks (each with a distinguished marked point in its interior) along each boundary component.  One records the following {\em combination data}:

\bi
\item the mapping tree $\mathbf{T}$, recording the deployment of the elements of $\AAA_0$ in the marked sphere $(S^2, Q)$; it is equipped with a self-map;
\item $\FFF^\SSS: \SSS_1 \to \SSS_0$, an induced map of a disjoint collection of spheres; 
\item $\FFF^\AAA: \AAA_1 \to \AAA_0$, an induced map of a disjoint collection of annuli;
\item set-theoretic gluing data $\tau$ needed to reconstruct the dynamics of $f$ on $Q$ from that of $\FFF^\SSS$ on the set $Q \intersect (S^2 \setminus \AAA_0) \union \ZZZ$ where $\ZZZ$ is the collection of distinguished added points, one for each disk; 
\item topological gluing data describing how to recover the original sphere from $\SSS_1$ and $\AAA_1$ and $\SSS_0$ and $\AAA_0$.
\ib

The Decomposition Theorem \cite[Thm. 5.1]{kmp:cds} asserts that up to combinatorial equivalence, $f$ can be reconstructed from the above combination data.  

There are notions of combinatorial equivalence for families of sphere and for annulus maps $\FFF^\SSS, \FFF^\AAA$.  

If $f$ is rational, the family $\FFF^\SSS$ is realized by a family of rational maps.  Their degrees and the size of their postcritical set are bounded in terms of $d$ and $q$.  By \cite[Cor. 3.7]{kmp:cds} as $\Gamma$ varies, up to combinatorial equivalence there are only finitely many such families obtained from $f$ by decomposition (the possibility of $\FFF^\SSS$ containing a cycle corresponding to a flexible Latt\`es example is not yet excluded; see below).  

If $f$ is rational with hyperbolic orbifold, the Thurston linear transformation $f_\Gamma$ does not have $1$ as an eigenvalue.  By \cite[Lemma 8.5]{kmp:cds}, this implies that as $\Gamma$ varies, up to combinatorial equivalence, there are only finitely many such families of annulus maps $\FFF^\AAA$.   

Suppose now $f$ has degree $d$ and $\#Q = q$.    The number of mapping trees $\mathbf{T}$ is bounded by a constant involving only $d$ and $q$.  The number of set-theoretic gluing maps  $\tau$ is similarly bounded. 

The preceding three paragraphs show that if $f$ is rational, then upon decomposing along invariant multicurves, up to combinatorial equivalence, the data needed to define $f$ as a combination range over finite sets.  

Suppose now $\Gamma_1, \Gamma_2$ are two invariant multicurves such that decomposing along $\Gamma_1$ and $\Gamma_2$ yield isomorphic mapping trees, set-theoretic and topological gluing data, and combinatorially equivalent sphere and annulus maps.  By the uniqueness of combinations theorem \cite[Thm. 4.5]{kmp:cds}, a pair $h^\SSS, h^\AAA$ of combinatorial equivalences between sphere and annulus maps, respectively, yields, upon combining, a combinatorial equivalence between the new, glued maps, and hence a combinatorial equivalence $h$ from $f$ to itself which represents an element of the pure mapping class group $\Mod(S^2, Q)$ and which sends $\Gamma_1$ to $\Gamma_2$.  

We now argue that $h$ is trivial and hence that $\Gamma_1 = \Gamma_2$.  
By the naturality of Thurston's pullback map used in the proof of the characterization theorem \cite[Prop. 2.1]{DH1}, the equivalence $h$ conjugates the pullback map $\sigma_f$ on Teichm\"uller space modelled on $(S^2, Q)$ to itself.  Since $f$ is rational with hyperbolic orbifold, $\sigma_f$ has a unique fixed-point $\tau$ which is therefore fixed by $h$, i.e. in the classification of mapping classes, $h$ is elliptic.  By Theorem \ref{thm:Thurston classification}, $h=\id$.  
\qed

\section{Analysis of the Rabbit}

In this and the next two sections, we use the spirit of the proof of Theorem \ref{thm:contracting implies eventually cycles} to analyze the pullback relation for quadratic polynomials with three finite postcritical points.  In two cases, the virtual endomorphism is contracting.  In the case of the dendrite $f(z)=z^2+i$, however, $\phi_f$ is not contracting on the correponding mapping class group $G$.  Nevertheless, the conclusion of Theorem \ref{thm:contracting implies eventually cycles} still holds. 
\gap

\noindent{\bf Remark.}  
A polynomial of the form $f(z)=z^d+c$ for which $\#(P_f \setminus \{\infty\}) = 3$  has no invariant (multi) curves.  Such a curve would bound a disk $D$ which is totally invariant up to isotopy relative to $P_f$.  This in turn would imply  that $P_f \setminus \{\infty\}$ would contain two points---one inside $D$, one outside---in distinct grand orbits.  This is impossible since $f$ has a single finite critical point.  
\gap
 
To keep the present notation close to that of \cite{bartholdi:nekrashevych:twisted}, we use $\psi$, not $\phi$, to denote virtual endomorphisms on the mapping class group.  The situation simplifies: since $\#P_f=4$ in each case, we have $\CCC = \MMM\CCC$, so that the pullback relation and pullback function coincide.

In this section, we prove Theorem \ref{thm:nwc_for_rabbit}.  Let $f(z)=z^2+c$ where $c$ is the unique complex parameter for which the origin is periodic of period $3$ and $\Imag(c)>0$.  Let $\CCC = \CCC(\rs, P_f)$ and $G=\Mod(\rs, P_f)$.  

Let $x$ and $y$ denote respectively the right Dehn twists about  the curves labelled $S$ and $T$ in Figure 2 of \cite{bartholdi:nekrashevych:twisted}.  For convenience, set $z=x^{-1}y^{-1}$; recalling the notational conventions, this means a left Dehn twist about $S$ is performed first, followed by a left Dehn twist about $T$ .  Then $z$ is a right Dehn twist about a curve separating $\{c,0\}$ from $\{c^2+c, \infty\}$.  

Any nontrivial element of $\Tw(\rs, P_f)$ is then uniquely expressible in one of the three forms 
\[ (x^n)^w, (y^n)^w, (z^n)^w\]
since the core curve of any twist can be mapped via a mapping class element $w$ to the core curve of either $x$, $y$, or $z$. 

The induced virtual endomorphism on $G$ is calculated in \cite{bartholdi:nekrashevych:twisted}; cf. also \cite[\S 6.6]{nekrashevych:book:selfsimilar}.  It is given by 
\[ \psi(x)=y, \;\;\;\psi(y^2)=x^{-1}y^{-1}=z, \;\;\; \psi(x^y)=1;\]
recall $x^y=y^{-1}xy$, so conjugation acts as a right action.   Let $H=\dom(\psi)=\genby{x, y^2, z}$.  
Note that 
\[ x^4 \mapsto y^4 \mapsto z^2 \mapsto x\]
and thus by Theorem \ref{thm:phi_and_LLL} under backwards iteration the core curves of $x, y, z$ form  a three-cycle.   

For convenience, set $u=yxy^{-1}$.  

\begin{lemma}
\label{lemma:3flavors}
We have 
\[
\begin{array}{lcl}
\psi((x^n)^w) & = & 
\left\{
\begin{array}{ll}
(y^n)^{\psi(w)}, & w \in H\\
1, & w \not\in H
\end{array}
\right.
\\
\psi((y^{2n})^w) & = & 
\left\{
\begin{array}{ll}
(z^n)^{\psi(w)}, & w \in H\\
(z^n)^{\psi(y^{-1}w)}, & w \not\in H
\end{array}
\right.
\\
\psi((z^{2n})^w) & = & 
\left\{
\begin{array}{ll}
(x^n)^{\psi(w)}, & w \in H\\
u^{\psi(y^{-1}w)}, & w \not\in H
\end{array}
\right.
\\
\psi((u^n)^w) & = &
\left\{
\begin{array}{ll}
1, & w \in H\\
(y^n)^{\psi(y^{-1}w)}, & w \not\in H
\end{array}
\right.
\end{array}
\]

\end{lemma}

\pf If $w \not\in H$ then $w=yw'$ where $w'\in H$, and 
\[
\begin{array}{ll}
 \psi((y^2)^w)=&\psi(w'^{-1}y^{-1}y^2yw')=\psi(w'^{-1}y^2w')=z^{\psi(w')}=z^{\psi(y^{-1}w)}\\
\psi((x^2)^w)=&\psi(w'^{-1}(y^{-1}xy)^2w')=\psi(w')^{-1}\cdot 1 \cdot \psi(w')=1\\
 \psi((z^2)^w)=&\psi(w'^{-1}y^{-1} z^2yw')=\psi(w'^{-1}y^{-1}\cdot x^{-1}y^{-1}x^{-1}y^{-1} \cdot y w')=\\
& \psi( w'^{-1} \cdot y^{-1}x^{-1}y \cdot y^{-2}\cdot x^{-1} \cdot w')=\\
& \psi(w')^{-1}\cdot 1 \cdot  yx \cdot y^{-1} \cdot \psi(w')=u^{\psi(y^{-1}w)},\\
\psi(u^w) =& \psi(w'^{-1}y^{-1} \cdot yxy^{-1} \cdot yw') = \psi(w')^{-1}\cdot \psi(x)\cdot \psi(w') = y^{\psi(y^{-1}w)}.
\end{array}
\]
If $w \in H$ then the desired identities follow directly from the definitions.
\qed

Define $\psihat: G \to G$ by setting $\psihat(g)=g$ for $g \in H$ and by $\psihat(g)=y^{-1}g$ otherwise.  

Let $\Sigma= \{x, y, z, u, 1\} \subset G$ 
and let $\pi: \Sigma \times G \to \CCC \union \{o\}$ be given by 
\[ \pi(c, w) = \mbox{the class of the core curve of the twist }c^w, \;\; c \in \Sigma, \; w \in G\]
where by convention $\pi(1, g)= o$ for all $g$.  

Define 
\[ E: \Sigma \times \Mod(S^2, P_f) \to \Sigma \]
by setting $E(c,w)$ to be the value given by the following table:
\[
\begin{array}{c||ccccc}
& x & y & z & u & 1 \\
\hline\hline 
w \in H & y & z & x & o & o  \\
w \not\in H & o & z & u & y & o.
\end{array}
\]
Lemma \ref{lemma:3flavors} and Theorem \ref{thm:phi_and_LLL} imply that the function $\mathbf{f}: \Sigma \times \Mod(S^2, P_f) \to \Sigma \times \Mod(S^2, P_f)$ given by $\mathbf{f}(c,w) = (E(c,w), \psihat(w))$ covers the pullback relation in the sense that the diagram 
\[ 
\begin{array}{ccc}
\Sigma  \times G& \stackrel{\mathbf{f}}{\longrightarrow} & \Sigma \times G
\\
\pi \downarrow & \; & \downarrow \pi \\
\CCC\union\{o\}& \stackrel{f}{\longleftarrow} &\CCC \union \{o\} \\
\end{array}
\] 
commutes.  Hence, we can lift iteration of the pullback relation on curves to iteration of $\mathbf{f}$. 
Theorem \ref{thm:nwc_for_rabbit} follows immediately from the above observations, Theorem \ref{thm:phi_and_LLL}, and the following Lemma, whose proof occupies the remainder of this section. 

\begin{lemma}
\label{lemma:psibar_eventually_trivial}
For any $w \in G$, there exists $n \in \N$ with $\psihat^{\circ n}(w)=1$.
\end{lemma}

\pf Identify the symmetric group on two symbols with the cyclic group of order two, $\Z_2$; the unique nontrivial element will be denoted $\sigma$. \footnote{The symbol ``1'', when within angle brackets, refers to the identity element of $G$.  When zero and one are used as a subscript on a vertical bar, they denote, respectively, the first and second $G$-coordinate in the wreath product.}   We begin by recalling the wreath recursion for the rabbit:
\[ \Phi: G \to (G \times G) \rtimes \Z_2 \]
given by 
\[ \Phi(x) =\genby{y, 1}, \;\;\; \Phi(y)=\genby{1, x^{-1}y^{-1} }\sigma; \]
see \cite{bartholdi:nekrashevych:twisted}.  

The map $\psihat$ is related to the recursion $\Phi$ in the following way: 
\begin{equation}
\label{eqn:psihat}
\psihat(g)=
\left\{
\begin{array}{lr}
g|_0& , g \in H\\
yx \cdot (g|_1) & ,g \not\in H.
\end{array}
\right.
\end{equation}
To see this, suppose first that $h \in H$.  Both $g \mapsto g|_0$ and $g \mapsto \psihat(g)$ are homomorphisms on $H$, so to show equality it is enough to show that they agree on the generators.  This is easily verified using the definitions.  If $g \not\in H$, then using the first case in (\ref{eqn:psihat}) we have 
\[
\begin{array}{ccl}
\psihat(g) & = & \psi(y^{-1}g)\\
\; & = & (y^{-1}\cdot g)|_0\\
\; & =&  \left. \left[ \genby{yx, 1}\sigma \cdot \genby{g|_0, g|_1}\sigma\right] \right|_0\\
\; & = & yx \cdot (g|_1).
\end{array}
\]
More generally, we have 
\gap

\noindent{\bf Claim.}  {\em For all $g \in G$ and $n \in \N$, }
\[ \psihat^{\circ n}(g)\in \left\{ g|_v, yx \cdot (g|_v), y \cdot (g|_v) : v \in \{0,1\}^n\right\}.\]

\pf We prove this by induction on $n$, the base case being handled by (\ref{eqn:psihat}).  Given $n \geq 2$ and $g \in G$, write 
\[ \psihat^{\circ n}(g) = \psihat( \underbrace{\psihat^{\circ (n-1)}(g)}_{a}) = \psihat(a).\]
By the inductive hypothesis, there exists $v \in \{0,1\}^{n-1}$ with 
\[ a \in \{ g|_v, y \cdot (g|_v), yx \cdot (g|_v) \}.\]
The proof breaks then into $3 \times 2$ cases, depending on the form of $a$ as an element of the above set, and whether or not $g|_v$ belongs to $H$.  The computations are tedious but straightforward; here is a representative calculation.  Suppose $a=yx \cdot (g|_v)$ and $g|_v = k \not\in H$.  Write $\Phi(k)=\genby{k_0, k_1}\sigma$.  Then $yx \cdot k \in H$ and so $\psihat(a) = \psi(yx\cdot k)$ is given by  
\[
\begin{array}{ccl}
\psi(yx\cdot k) &=& \left[\genby{1, x^{-1}y^{-1}}\sigma \cdot \genby{y,1}\cdot \genby{k_0, k_1}\sigma\right]|_0\\
\;&  = & \left[\genby{1, x^{-1}y^{-1}}\genby{1,y}\genby{k_1, k_0}\right]|_0\\
\; & = & k_1 = k|_1 = (g|_v)|_1\\
\; & = & g|_{v1}.
\end{array}
\]
\qedspecial{Claim}

In \cite{bartholdi:nekrashevych:twisted} it is shown that the recursion for the rabbit is contracting in the following sense:  there is a finite set $\NNN \subset G$ such that given any $g \in G$, there is an integer $N$ such that for any $n \geq N$, and any $v \in \{0,1\}^n$, we have $g|_v \in \NNN$.  They compute the set $\NNN$; it is given by 
\[ \NNN = \{1, x, y, yx, x^{-1}, y^{-1}, x^{-1}y^{-1}\}.\]
This and the preceding Claim imply that under iteration of $\psihat$, any $w \in G$ eventually lands in an element of the set 
\[ \NNN, \;\;y \cdot \NNN, \;\;yx \cdot \NNN.\]
The dynamics of $\psihat$ on this union is easily calculated by hand:
\[
\begin{array}{l}
yxyx \mapsto x^{-1} \mapsto y^{-1} \mapsto yx \mapsto y \mapsto 1\\
yx^2 \mapsto y^2 \mapsto x^{-1}y^{-1} \mapsto yx\\
yxy^{-1} \mapsto 1 \leftarrow yx^{-1}y^{-1}\\
y^2x \mapsto x^{-1}, yx^{-1}\mapsto y^{-1}, yxy \mapsto x^{-1}y^{-1}.
\end{array}
\]

\qed

\noindent{\bf Remarks.}  The argument given above is a direct imitation of the techniques in \cite{bartholdi:nekrashevych:twisted}, elaborating on their Proposition 4.2.  The relationship between the wreath recursion of the rabbit and our map $\psihat$ seems substantially different than that between the wreath recursion of the rabbit and their map $\hat{\psi}$: the map $\mathbf{f}$ formally defines an automaton over the alphabet $\Sigma$ with states $G$.  For a given $g \in G$, however, the corresponding map $\Sigma \to \Sigma$ need not be invertible.

\section{Analysis of $z^2+i$}

In this section, we prove Theorem \ref{thm:eventually_trivial_dendrite}.  Let $f(z)=z^2+i$ and let now $\CCC = \CCC(\rs, P_f)$.  Our analysis proceeds differently:   the virtual endomorphism $\psi$ and corresponding recursion is no longer contracting on $G=\Mod(\rs, P_f)$.   The analysis exploits the fact that the analysis of the action of iteration of $\psi$ on conjugacy classes can be reduced to the single case of iteration of $b^w \mapsto b^{\psihat(w)}$ where $\psihat$ is similarly defined.  To conclude the argument, however,  we exploit cancellations caused by identities such as $b^b = b$.  

The domain of the virtual endomorphism $\phi_f$ in this case is $H=\genby{a^2, b, b^a}$, and the action on the generators is given by 
\[ \psi(a^2)=1, \psi(b)=b^{-1}a^{-1}, \psi(b^a)=b.\] 
Let us set $c=b^{-1}a^{-1}$, so that the loop representing $c$ contains a peripheral disk bounding infinity to its left-hand side. 

\begin{lemma}
\label{lemma:3flavors_dendrite}
We have 
\[ \psi((a^{2n})^w)=1
\]
\[
\psi((c^{2n})^w)=
\left\{
\begin{array}{rl}
(a^n)^{b^{-1}\psi(w)}, & w \in H\\
(a^n)^{b^{-1}\psi(a^{-1}w)}, & w \not\in H
\end{array}
\right.
\]
\[
\psi((b^n)^w)=
\left\{
\begin{array}{rl}
(c^n)^{\psi(w)}, & w \in H\\
(b^n)^{\psi(a^{-1}w)}, & w \not\in H
\end{array}
\right.
\]
\end{lemma}

The proof is entirely analgous to that of Lemma \ref{lemma:3flavors}.  

Lemma \ref{lemma:3flavors_dendrite} implies that $\psi^{\circ 2}(u^w)$ is trivial whenever $u\in\{a,c\}$, and so $\psi^{\circ 3}(b^w)$ is trivial if $w \in \dom(\psi)$.  The long-term behavior of iterates of $\psi$ on elements of the form $b^w$ is therefore again dictated by iteration of the map $w \mapsto \psihat(w)$, where $\psihat(w)=w$, $w \in H$ and $\psihat(w)=\psi(a^{-1}w)$, $w \not\in H$.   

Let $|w|$ denote the word length of $w$ with respect to the generating set $\{a^{\pm 1}, b^{\pm 1}, c^{\pm 1}\}$.

\begin{lemma}
\label{lemma:nwc_for_dendrite_recursion}
The map $\psihat$ satisfies the following recursion relations:
\begin{equation}
\psihat(aw)=
\left\{
\begin{array}{rl}
1 \cdot \psihat(w), & w \in H \\
1 \cdot \psihat(w), & w \not\in H.
\end{array}
\right.
\end{equation}

\begin{equation}
\psihat(a^{-1}w)=
\left\{
\begin{array}{rl}
1\cdot \psihat(w), & w \in H\\
1\cdot \psihat(w), & w \not\in H.
\end{array}
\right.
\end{equation}
\begin{equation}
\psihat(b^{\pm 1}w)=
\left\{
\begin{array}{rl}
c^{\pm 1} \cdot \psihat(w), & w \in H\\
b^{\pm 1} \cdot \psihat(w), & w \not\in H.
\end{array}
\right.
\end{equation}
\begin{equation}
\psihat(cw)=
\left\{
\begin{array}{rl}
b^{-1}\cdot \psihat(w), & w \in H\\
c^{-1}\cdot \psihat(w), & w \not\in H.
\end{array}
\right.
\end{equation}
\begin{equation}
\psihat(c^{-1}w)=
\left\{
\begin{array}{rl}
c \cdot \psihat(w), & w \in H\\
b \cdot \psihat(w), & w \not\in H.
\end{array}
\right.
\end{equation}
In particular, $|\psihat(w)| \leq |w|$ for all $w$.
\end{lemma}

The computations are straightforward.

\begin{lemma}
\label{lemma:dendrite_contracting}
If $w\not\in H$ and $\psihat(w) \not\in H$, then there exists $v$ such that 
$\psi^{\circ 2}(b^w)=b^v$ and $|v|<|w|$.
\end{lemma}

\pf Let $w=s_1s_2\ldots s_l$ be a minimal length representation of $w$ as a word in the generators $\{a^{\pm 1}, b^{\pm 1}, c^{\pm 1}\}$.
If $s_1 = b^{\pm 1}$ then set $w'=s_2s_3 \ldots  s_l$ and observe 
\[ \psi(b^w) = \psi(b^{b^{\pm 1}w'}) = \psi(b^{w'}) = b^{\psihat(w')}\]
and so 
\[ \psi^{\circ 2}(b^w)=b^v, \;\;\; v = \psihat^{\circ 2}(w').\]
By the previous Lemma 8.2, 
\[ |v| = |\psihat^{\circ 2}(w')| = |w'| \leq l-1 < |w|.\]

If $s_1 = a^{\pm 1}$ or $s_1 = c$ then a similar calculation shows the desired inequality.  

If $s_1 = c^{-1}$ we argue as follows.  Set $w=c^{-1}w'$ where $w'$ is as above.  Then 
\[ \psi(b^w) = \psi(b^{abw'}) =b^{\psi(bw')}=b^{c\psi(w')}=b^{a^{-1}\psi(w')}=b^{a^{-1}\psihat(w')}\]
and so 
\[ \psi^{\circ 2}(b^w)=\psi(b^{a^{-1}\psihat(w')})=b^{\psihat^{\circ 2}(w')}=b^v, \;\;\; v = \psihat^{\circ 2}(w').\]
Again $|v|=|\psihat^{\circ 2}(w')| \leq |w'| = l-1 < |w|$, as required.  
\qed

Now let $\gamma \in \CCC$ be an arbitrary curve.  Then the Dehn twist about $\gamma$ is given by $u^w$ where $u \in \{a, b, c\}$.  By Lemma \ref{lemma:3flavors_dendrite},  if $u \in \{a,c\}$ then $\psi^{\circ 2}((u^4)^w))=1$, hence all preimages of $\gamma$ under $f^{\circ 2}$ are trivial.  So we now assume $u=b$.   Let $k=2|w|$.   By induction and Lemmas \ref{lemma:3flavors_dendrite} and \ref{lemma:dendrite_contracting},  $\psi^{\circ i}((b^4)^w)$ is either trivial already, a conjugate of $a$ or of $c$ (in which case it becomes trivial upon applying $\psi^{\circ 2}$), or is equal to $b^4$ (in which case it becomes trivial upon applying $\psi^{\circ 3}$).  
\qed

\section{Preperiod 1, period 2}

In this section, we prove Theorem \ref{thm:eventually_trivial_pre1per2}.  The computations were carried out by R. Lodge.    

There is a unique complex number $c \approx -0.2282\ldots + 1.1151\ldots i $ with $\Imag(c)>0$ for which, under iteration of $f_{1/4}(z)=f(z):=z^2+c$, the origin after one iteration lands in a two-cycle; it is the landing point of the angle $1/4$ external ray of the Mandelbrot set.  

Let $a$ and $b$ denote respectively the right Dehn twists about the curves labelled ``$a$'' and ``$b$'' in Figure 18 of \cite{bartholdi:nekrashevych:twisted}.  The domain $H$ of the virtual endomorphism $\psi$ is given by $H=\genby{a^2, b, aba^{-1}=b^{a^{-1}} }$, and the action on the generators is given by 
\[ \psi(a^2)=b, \;\; \psi(b)=b^{-1}a^{-1}, \;\; \psi(aba^{-1})=a.\]
Define $\psibar: G \to G$ by 
\[ \psibar(w) = \left\{
\begin{array}{rr}
\psi(w), & w \in H\\
\psi(a^{-1}w), & w \in aH.
\end{array}
\right.
\]
It is convenient to set $c=b^{-1}a^{-1}, d=b^{-1}ab, e=aba^{-1}, f=b^{-1}a^{-1}bab$.  
As before, the pullback relation on curves lifts to a function $\mathbf{f}: \Sigma \times G \to \Sigma \times G$. If $g \in \Sigma :=\{1, a, b, c, d, e, f\}$ the smallest positive integer $k$ for which $g^k \in H$ is given by $k=1, 2, 1, 2, 2, 1, 1$, respectively, as can be seen by counting (with sign) the powers of the generator $a$ appearing in an expression for $g$.   The corresponding table summarizing the values of the corresponding function $E$ in this case is given by 
\[
\begin{array}{c||cccccc}
& a & b & c & d & e & f  \\
\hline\hline 
w \in H & 		b & c & 1 & e & a & a  \\
w \not\in H & 	b & d & 1 & f & c & d.
\end{array}
\]
By computations entirely analogous to those in \S 8, we find that for all $g \in G$, 
\[ \psibar(g) = \left\{
\begin{array}{rr}
g|_0, & w \in H\\
b^{-1}(g|_1), & w \in aH
\end{array}
\right.
\]
where $\Phi(g)=\genby{g|_0, g|_1}\sigma$ is the wreath recursion given on generators by 
\[ \Phi(a) = \genby{1, b}\sigma, \;\; \Phi(b)=\genby{b^{-1}a^{-1}, a}.\]
By induction, one finds that for all $g \in G$ and $n \in \N$ that 
\[ \psibar^{\circ n}(g) \in \left\{ g_v, b^{-1}(g|_v), ab(g|_v), b^{-1}a^{-1}(g|_v), a(g|_v) : v \in \{0,1\}^n\right\};\]
the proof is a tedious but straightforward computation involving ten cases depending the coset containing $g$ and the form of $\psibar^{\circ (n-1)}(g)$.  From \cite{bartholdi:nekrashevych:twisted}, the nucleus of the wreath recursion is given by 
\[ \NNN = \left\{1, a^{\pm 1}, b^{\pm 1}, (ab)^{\pm 1}, (a^{-1}b)^{\pm}\}\right\} \]
and so under iteration of $\psibar$ every element eventually lands in the union of one of the five sets 
\[ \NNN, \;b^{-1}\NNN, \;ab\NNN, \;b^{-1}a^{-1}\NNN, \;a\NNN.\]
Another straightforward computation shows that if $n \geq 6$ then every element $g$ lying in one of these five sets satisfies $\psibar^{\circ n}(g) \in \{1, b^{-1}a^{-1}\}$.  

It follows that under iterated pullback, every curve orbit  is covered by an orbit of $\mathbf{f}$ that lands in the set $\Sigma \union \Sigma^{b^{-1}a^{-1}}$.  Another straightforward computation shows that upon pulling back five times, such a curve becomes trivial.
\qed

\def\cprime{$'$}

\me

\end{document}